\documentclass[12pt]{article}

\usepackage{amsmath}
\usepackage{amssymb}
\usepackage{amsthm}
\usepackage{amscd}
\usepackage{amsfonts}
\usepackage{dsfont}
\usepackage{cases}
\usepackage[applemac]{inputenc}
\usepackage{enumerate}
\usepackage[all]{xy}
\usepackage{esint}
\usepackage{graphicx}
\usepackage{url}
\theoremstyle{plain}

\theoremstyle{definition}

\newcommand{\N}{\mathbb N}

\newcommand{\E}{\mathrm E}

\date{}
\begin{document}

\title{Odds -Theorem and Monotonicity}
\author{F. Thomas Bruss \\Universit\'e Libre de Bruxelles}
\maketitle
\begin{abstract} \noindent Given a finite sequence of events and a well-defined notion of events being {\it interesting}, the Odds-theorem  (Bruss (2000)) gives an online strategy to stop on the last interesting event. It is optimal for independent events.  Here we study questions in how far optimal win probabilities mirror monotonicity properties of the underlying sequence of probabilities of events. We make these questions precise, motivate them, and then give complete answers. This note, concentrating on the original Odds-theorem, is elementary, and the answers are hoped to be of interest. We include several applications.  \noindent\medskip

\medskip
\noindent {\bf Keywords}: Odds-algorithm; Secretary problem; Selection criteria; Multiple stopping problems; Group interviews; Games; Clinical trial; Prophet inequality; 

\bigskip
\noindent {\bf Math. Subj. classif.} 60G40
\end{abstract}
\section {Motivation}
The Odds-theorem is an easy-to-apply theorem in the domain of optimal stopping. It gives an online strategy to stop on the very last {\it interesting} event
of a given sequence of events. Its interest lies in the flexibility of the notion {\it interesting},  in its optimality for independent events, and in the {\it odds-algorithm} which flows from it. The Odds-theorem can also be useful for a similar setting with conditionally independent events. If not stated otherwise, we assume independence. In this case the odds-algorithm is also optimal as an algorithm. 

The present article is a second addendum to the Odds-Theorem, after the short note of B. (2003) proving the general lower bound $1/e$ for the win probability. ("B." stands throughout all references for the author's name.) In this article we mainly examine questions of monotonicity. Monotonicity questions arise, among other instances, when the decision maker has some influence on the order of events within the sequence, or when the length of the sequence may vary after observations have begun.  \subsection{Examples and related work}
Questions in how far optimal win probabilities mirror monotonicity properties of the underlying sequence of probabilities of events are shown to be relevant in several contexts. 
As we will see in the Section Applications
(Section 5), they can help us to (quickly) decide which game we want to play if we have a choice,  or may give us advice in scheduling interviews of candidates in order to make a better choice. They may  help us to find  answers  for new questions concerning well-known selection problems,  but also be equally important in quite different contexts, and also  from a different point of view, as for instance in planning or re-organizing clinical trials or specific arrangements for sequential medical treatments.

We will focus for the major part of the present article on the original Odds-theorem, including its continuous-time version for processes with independent increments presented in the same paper (B. (2000)). We should add that now there exist  several interesting variations of the underlying model and/or its payoff-function. These  include multiple stopping problems, such as e.g. in Matsui and Ano (2014), (2016), or modified payoffs (see e.g. Tamaki (2010), (2011)), or  again modified in such a way that they may be helpful for related game  problems in continuous-time, such as in Szajowski (2007). (For a best-choice problem with dependent criteria see e.g. Gnedin (1994).) Moreover, the Odds-Theorem can be adapted for conditional independent events, as in Ferguson (2016) and other papers. Such models involve statistical learning, in which monotonicity properties
may also play a role. For a review of results on developments of the Odds Theorem see Dendievel (2013).

We also want to add that an interesting alternative proof of the Odds-theorem  has recently be given in  Goldenshluger et al. (2019). 

\section{From $n$ fixed to $n$ varying.} We recall the notation of B. (2000).
 Let $(E_k)_{k=1, 2, \cdots ,n}$ be a sequence of $n$ events defined on some probability space $(\Omega,{\cal A}, P)$. Suppose we have a well-defined criterion according to which $E_k$ is an {\it interesting} event. Let $I_k={\bf 1}\{E_k~{\rm is~interesting}\}.$ Correspondingly, we call   $P(I_k=1)/P(I_k=0)$ the {\it odds} of $E_k$ being interesting. Further, let $p_k=P(I_k=1), q_k=1-p_k, r_k=p_k/q_k, ~ k = 1, 2, \cdots, n,$ and let \begin{align}Q_k=\prod_{j=k}^n q_j;~ R_k=\sum_{j=k}^nr_j.\end{align}
Then, if the $I_k$ are independent, the index $k$  maximizing $Q_kR_k$, denoted by $s$, yields the optimal value $V(n)$, namely $V(n)= Q_sR_s.$ Here $s=1$ if $R_1<1,$ otherwise $s$ is the largest $k$ such that
$R_k\ge 1.$ This is the Odds-theorem, and this $s$ is called {\it optimal threshold. }

\medskip \noindent It is again useful to point out that the answer $V(n)=Q_sR_s$ is complete. It also covers  the case $Q_s=0$ and $R_s=\infty$ since $V(n)$ stays well-defined, as shown below.

\bigskip
\noindent{\bf Corollary 2.1} If, for $n$ fixed,  $s$ is the optimal threshold, and $Q_s=0,$ then $V(n)=Q_{s+1}.$

\medskip
\noindent{\bf Proof:~}
If $s$ is optimal with $p_s=1$, then clearly $s$ must be  the {\it last} $j$ such that $p_j=1.$ With $1=p_s=r_s\, q_s$ the  product  $q_s\,r_s$ is no undetermined form. Therefore from (1) and (2), using $R_s=R_{s+1}+r_s, ~Q_s=q_sQ_{s+1},$ 
\begin{align}Q_s R_s=Q_{s+1} \,(q_sR_{s+1}+q_sr_s)=Q_{s+1}(q_sR_{s+1}+1)=Q_{s+1}=\prod_{j=s+1}^nq_j,\end{align} because $R_{s+1}<1$ by definition of $s$, and $q_s=1-p_s=0.$
The answer stays also correct for $s=n$ with the standard definition that an empty product equals $1.$ \qed

\medskip \noindent {\bf Remark 2.1} We conclude from Corollary 2.1 that, whenever we deal with $Q_sR_s,$ we can always assume $Q_s>0$; otherwise $Q_sR_s$ reduces to $Q_{s+1}>0.$ This simple result will be used repeatedly in what follows.

\subsection {The $n$-problem}
In the following we define the setting of the Odds-theorem for $n$ varying. Although it seems evident what is meant by $n$ varying, a clear terminology will  keep our arguments simple.

We speak of a {\it $n$-problem} for an {\it underlying} sequence of probabilities $p_1, p_2, \cdots,$ if the problem of stopping on the last interesting event applies to the stopped sequence $I_1, \cdots, I_n$ with $\E(I_k)=p_k, k=1, 2, \cdots, n.$ More precisely:

\medskip
(i)
We say that we {\it win} in the $n$-{\it problem} if we succeed to stop on the last index $k\in \{1,2,\cdots ,n\}$ with $I_k=1.$ 

\medskip
(ii) A stopping time $\sigma$ is said to be {\it optimal} for the $n$-problem if $\sigma$ maximizes the win probability for the $n$-problem. The corresponding value is denoted by $V(n).$

\medskip
(iii) We say that $s(n)$ is the optimal threshold for the $n$-problem, if the stopping time
\begin{align}\sigma_n= \inf\{s(n) \le k \le n: I_k=1\}\wedge n\end{align}solves the $n$-problem. Here, as for $n$ fixed, it is understood that one cannot return to an $I_j=1$ passed over before, and that the stopping time is also non-anticipative, i.e. $\{\sigma_n = k\} \in {\cal F}_k$ where ${\cal F}_k$ denotes the $\sigma$-field generated by $I_1, I_2, \cdots,I_k$. 

\medskip
From the Odds-theorem (B.(2000) p. 1385) we have therefore \begin{align}V(n)=Q(s,n)R(s,n),\end{align}where, 
for $1\le k \le n$,
\begin{align}Q(k,n)=\prod_{j=k}^n q_j; ~R(k,n)
=\sum_{j=k}^n r_j,\end{align} 
and \begin{align}s:=s(n)=
\begin{cases}1&,{\rm~ if}~ R(1,n)\le 1\\
\sup\{1\le k\le n:R(k,n) \ge 1\}&, \mbox{~otherwise}.
\end{cases}\end{align}

\medskip \noindent According to (6), we will use the simplified notation $Q(s,n):=Q(s(n),n) $ and $R(s,n):=R(s(n),n)$ whenever this is not ambiguous.

\section{Monotonicity}

We are now ready to tackle questions of interest concerning the monotonicity of $V(n)$ in $n.$  Since we have an explicit and simple formula for $V(n),$ such questions, including when $V(n)=V(n+j)$ for fixed $j\in \N$ will hold, are not deep, of course. Our focus will be on trying to see certain facts quickly, and what their implications are. Also, our objective is to increase intuition of what will happen if the underlying sequence changes in a certain way.

\medskip The following Lemma 3.1 and Theorem 3.1
are the basic result.

\bigskip
\noindent{\bf Lemma 3.1} For an underlying sequence  $p_1, p_2, ....$ with $0\le p_j\le 1$ for all $j,$ let\begin{align} N^*=\sup\{n \in \N: R(1,n) \le 1\}.\end{align} Then the optimal win probability $V(n)$ is non-decreasing in $n$ for $1\le n \le N^*.$

\bigskip
\noindent{\bf Proof.}~ Since $s(n)=1$ for $n\in\{1, 2, \cdots, N^*\}$ by definition of $N^*$, we have from the optimality of the threshold $s(n)$ the value $V(n)=Q(1,n)\,R(1,n)$ for $n\le N^*.$ Thus, by definition of $Q(k,n)$ and $R(k,n)$ in (1) and (2), we get  for $1\le n< N^*$
\begin{align*}V(n+1)=Q(1, n+1)\,R(1,n+1)=Q(1,n)\,q_{n+1}\,\big(R(1,n)+r_{n+1}\big)\\
=Q(1,n)\,\big(q_{n+1}R(1,n)+p_{n+1}\big)\ge Q(1,n)\,R(1,n)\,(q_{n+1}+p_{n+1})\\= Q(1,n)\,R(1,n)=V(n)\end{align*} where the inequality follows from $R(1,n)\le1.$ Thus $V(n+1)\ge V(n)$ for $n<N^*.$\qed

\bigskip
The index $N^*$ defined in Lemma 1 is a benchmark in the sense
that assumptions for $p_1, p_2, \cdots, p_{N^*}$ are irrelevant for the monotone behaviour of $V(n).$ As we shall see in the following, from $N^*$ onward, $V(n)$  mimicks
(simple) monotonicity assumptions of the $p_n$ on $\{N^*, N^*+1, \cdots\}.$ We will see later on that the latter is not necessarily true for {\it strict} monotonicity.
\bigskip

\noindent {\bf Theorem 3.1} The sequence of optimal values $V(n)_{n \ge N^*}$ for the $n$-problems reflects monotonicity properties of an underlying sequence $(p_n)_{n \ge N^*}$ as follows:

\medskip
\medskip
(A) If the sequence $(p_n)_{n\ge N^*}$ is non-increasing then the optimal values are non-increasing for $n\ge N^*.$

\medskip
(B) If $(p_n)_{n\ge N^*}$ is non-decreasing then $V(n)$ is non-decreasing for all $n\in \N$.

\bigskip

\noindent 
{\bf Proof.}
\medskip

\noindent{\bf (A)}  ~~We first show that if the success probabilities $p_j$  are non-increasing then the optimal threshold $s(n)$ for the $n$-problem satisfies\begin{align} \forall \,n\in\{1, 2, \cdots\}: s(n)=s \implies s(n+1)\in \{s, s+1\}.\end{align}

Indeed, we first note that by definition $1\le s(j)\le j$ and  $s(j)\le s(j+1)$ since all odds are non-negative. Also, $R(s,n)-1< r_s$ since, from (6), $R(s+1,n)=R(s,n)-r_s<1$ and $R(s,n)\ge 1$. Moreover, with non-increasing $p_j$ we see that we have non-increasing odds $r_j=p_j/q_j.$ Consequently, there are only two possibilities by passing from $n$ to $n+1$: if $R(s+1,n+1)=R(s,n)-r_s+r_{n+1}\ge 1$ we get $s(n+1)=s(n)+1=s+1$, otherwise $s(n+1)=s.$  This proves statement (8).
\smallskip

\noindent Hence we  have to consider for the proof of (A) only two cases, namely (i) $s(n+1)=s(n),$ and (ii) $s(n+1)=s(n)+1.$  Clearly, we can limit our interest to $n\ge N^*.$

\bigskip 

(i) Let  $s(n+1)=s.$ Then (4) and (5) imply that the inequality $V(n+1)\le V(n)$ is equivalent to the inequality\begin{align} Q(s,n+1)R(s,n+1)\le Q(s,n)R(s,n),\end{align}which, according to (5), is again equivalent to 
\begin{align}q_{n+1}Q(s,n)\,(R(s,n)+r_{n+1}) \le Q(s,n)R(s,n).\end{align}
Recall Remark 2.1 and divide by $Q(s,n)>0.$ Using  $ r_{n+1}=p_{n+1}/q_{n+1}$ we see  that inequality  (10) becomes \begin{align}q_{n+1}R(s,n))+p_{n+1}\le R(s,n).\end{align} This inequality is always true for $n\ge N^*$, since  $p_k+q_k=1$ for all $k,$ and since for $n\ge N^*$ we have $R(s,n)\ge 1$ by definition of $s:=s(n).$

\bigskip\medskip (ii)
If $s(n+1)=s+1$ then we must prove that the condition  \begin{align}Q(s+1,n+1)R(s+1,n+1)\le Q(s,n)R(s,n)\end{align} will hold for $n\ge N^*. $ This is slightly more involved.

\smallskip \noindent
By definition of $Q(s,n)$ and $R(n,s)$ the inequality (12) is now equivalent to
\begin{align} \frac{q_{n+1}}{q_s}Q(s,n)\,(R(s,n)+r_{n+1}-r_s) \le Q(s,n)R(s,n).\end{align} We first note that the case $q_s=q_{n+1}$ is trivial because then we have also $r_{n+1}=r_s$ so that both sides of (13) become $Q(s,n)R(s,n),$ and thus the statement is true.

\smallskip
Hence we can confine our interest to $q_s\not=q_{n+1}.$ Since we assumed the $p_j$ non-increasing, this means  $p_s>p_{n+1}$ and $r_s>r_{n+1}.$

Independently, we have seen already that we can focus our interest on $Q(n,s)>0$, implying $0<p_s<1$ and  $0<q_{n+1}<1.$ 
Therefore, dividing inequality (13) by $Q(s,n) > 0 $ and multiplying it by $q_{s}>0,$ it becomes
\begin{align}R(s,n)\left(q_s-q_{n+1}\right)\ge q_{n+1}(r_{n+1}-r_s).\end{align} Since the rhs of (14) can be written $p_{n+1}
-q_{n+1} (p_s/q_s)$ we obtain, using $q_s-q_{n+1}=p_{n+1}-p_s< 0,$  \begin{align}R(s,n)\le
\frac {-p_{n+1}+r_s q_{n+1}}{-p_{n+1}+p_s}.\end{align} With non-increasing  $p_j$'s we have non-decreasing $q_j$'s so that $$r_sq_{n+1}\ge r_sq_s=p_s.$$ Therefore the rhs of (15) is greater or equal to 1 as it should be in the non-trivial case $n\ge N^*$ by definition of $R(s,n)$.

\smallskip
However, here we have to observe an additional combined constraint. By passing from $n$ to $n+1,$
the optimal threshold index $s(n+1)$ for the $(n+1)$-problem becomes $s(n)+1$ if and only if
\begin{align} R(s(n),n) \ge 1, {~\rm and~} R(s(n)+1,n) < 1,  {~\rm and~} R(s(n)+1,n+1)\ge 1.\end{align}
Since $R(s(n),n)-R(s(n)+1,n)=r_{s(n)}$ and 
$R(s(n)+1,n+1)-R(s(n)+1,n)=r_{n+1}\ge r_{s(n)}$
the constraints in (16) are satisfied  if the rhs of inequality (15) does not exceed or reach the value $1+r_{s(n)}.$ Indeed, we see that  inequality (15) is sharp, namely,
with $s:=s(n)$,
\begin{align} \frac{r_sq_{n+1}-p_{n+1}}{p_s-p_{n+1}}=1+r_s.\end{align} To see this, note that $1+r_s$ can be written as $q_s^{-1}.$ Since $q_s>0$, the equation (17) is equivalent to  $q_{n+1} p_s-p_{n+1}q_s=p_s-p_{n+1},$ and this is straightforward to verify. 

\smallskip
This completes the proof of part (A).

\bigskip \noindent {\bf (B)}~~Suppose now that the sequence $(p_n)_{n\ge N_0}$ is non-decreasing. Note that, although this means that the $q_n$ are non-increasing, we cannot use here a duality argument based on re-interpreting the $I_k$ as $1-I_k.$ Therefore, the proof of (B) does not follow directly from the proof of (A), but we can use several parts of it.

\smallskip
We first note that the $(n+1)$-optimal threshold $s(n+1)$ can now no longer coincide with $s(n).$ In fact, with $p_{n+1}\ge p_{s(n)}, $ and thus 
$r_{n+1}\ge r_{s(n)},$ we see from (6) that the sum of odds $R(s(n+1), n+1)$ would otherwise not be the minimal tail sum of odds to reach or exceed $1.$ This implies that the part  (i) of the proof of (A)
is now irrelevant, and that (ii) should now read $s(n+1)\ge s(n)+1.$ 

\medskip
To begin with $s(n+1)\ge s(n)+1$, suppose first that $s(n+1)=s(n)+1.$ 
Then we can use the proof of the part (A)  literally by reversing all inequality signs in the equivalence (9)-(13). Also, the equality (17) stays  valid. Hence 
\begin{align}V(n+1)=Q(s(n)+1,n+1)R(s(n)+1,n+1)\ge Q(s(n),n)R(s(n),n)=V(n),\end{align} so that the statement (B) is proved for $s(n+1)=s(n)+1.$ 

\smallskip
Furthermore, we can now use an important part of the proof of the Odds-theorem B.(2000). It is the part dealing with uni-modality (see p. 1386, line 3 up to equation (4)). 
 It was shown there that
the function $f(k,n):= Q(k,n)R(k, n)$ is, for fixed $n$, unimodal in $k.$ The uni-modality holds in the sense that $f(k,n)$ is either non-increasing for all $1\le k \le n$, or else non-decreasing up to its maximum, and non-increasing thereafter. 

Now, replacing $n$ by $n+1$ we know then from the inequality figuring in (18) that the index $k=s(n)+1$ must belong to the non-decreasing wing of the graph of $f(k,n+1)$. Since $s(n+1)=\arg \max_{1\le k\le s(n+1} \{f(k,n+1)\}$, this  uni-modality property implies \begin{align}f(s(n)+1,n+1)\le f(s(n)+2,n+1)\le\cdots\le f(s(n+1), n+1).\end{align}
But then, given that the rhs of (19) corresponds to $f(s(n+1), n+1)=V(n+1)$  (by definition of $s(n+1)$) whereas  the inequality(18) affirms that $V(n) \le f(s(n)+1,n+1).$  This is true since the index $s(n)+1$ with $s(n)\le s(n)+1\le s(n+1)$ lies in the non-decreasing part  of $f(\cdot, n+1)).$ Hence, taking both together we have $$V(n)\le f(s(n)+1,n+1)\le V(n+1).$$

\medskip \noindent This proves part (B) and completes the proof of Theorem 3.1.\qed

\bigskip
\noindent The following easy observation is worth pointing out.

\bigskip
\noindent {\bf Corollary 3.1} $\forall n\in\{1, 2, \cdots\}: s(n+1)=s(n) \implies V(n+1)\le V(n)$

\medskip
\noindent {\bf Proof.} In the part A (i) of Theorem 3.1
we only used $R(s(n),n)\ge 1$ for $n\ge N^*.$ However, the latter holds by definition and is independent of  monotonicity assumptions (although the hypothesis $s(n+1)=s(n)$ itself is not, as just seen before). \qed
\section{Uniqueness of optimal thresholds and values}
Corollary 1 of B.\,(2000) (see p. 1387) says for fixed $n$: $V(n)=Q_sR_s=Q_{s-1}R_{s-1}$ if and only if $R_s=1.$ This translates for the $n$-problem based on an underlying sequence $(p_j)$ directly into:

\medskip
\noindent{\bf Corollary 4.1} \begin{align}V(n)=Q(s(n)-1,n)R(s(n)-1,n) \iff 
R(s(n),n)=1,\end{align} that is, two consecutive indices $s-1$ and $s$ are both optimal thresholds for the $n$-problem if and only if the sum of the relevant odds in the  $n$-problem equals $1.$

\bigskip \noindent The following results are complementary to the preceding one. We give criteria for values of different $n$-problems to coincide. Since we have an explicit formula  for $V(n)$ in terms of $p_1, p_2, \cdots,p_ n$ and $s(n),$ we have a straightforward equivalence, namely $V(n+1)=V(n)$ if and only if $Q(s(n+1), n+1) R(s(n+1), n+1)=Q(s(n),n)R(s(n),n).$  Since from (5), (by putting for $Q(a,b)=1$ and $R(a,b)=0$ for $b<a$),\begin{align}Q\big(s(n+1), n+1\big)=q_{n+1}Q\big(s(n),n\big)\,\Big/\prod_{j=s(n)}^{s(n+1)-1}q_j
\,=\,\frac{q_{n+1}Q\big(s(n),n\big)}{Q\big(s(n),s(n+1)-1\big)},\\
R\big(s(n+1),n+1\big)=R\big(s(n),n\big)-R\big(s(n),s(n+1)-1\big)+r_{n+1},~~~~~~~~~~~~~\end{align} we can solve the equation $V(n)=V(n+1)$ explicitly for $R(s(n),n)$. This requires to  compute $s(n+1),$ which is no problem, of course, but means additional work.
In the same way we could derive for an arbitrarily fixed $j\in\{1,2, \cdots\}$ a criterion for $V(n+j)
=V(n)$ to hold. It clearly suffices to adapt the above formulae in (21) and (22) and then to solve the equivalence equation again for  $R(s(n),n)$. 

\medskip As our primary goal is to increase the ease of application of the Odds-Theorem and to see implications as quickly as possible, (21) and (21) are slightly too complicated for this purpose. In the case of monotonicity things are simpler:

\bigskip
\noindent{\bf Theorem 4.1} If the underlying sequence $(p_j)$ with $0<p_j<1$  is non-increasing then\begin{align*}V(n+1)=V(n) \end{align*} if and only if one of the two following conditions are satisfied
$${\rm(a)}~R(s(n),n)=1{~\rm or~} R(s(n+1),n)=1~$$
$${\rm(b)}~p_{s(n)}=p_{n+1}.~~~~~~~~~~~~~~~~~~~~~~~~~~~~~~~~$$

\medskip\noindent{\bf Proof}:~We recall from the first part of the proof of Theorem 3.1  (see (8)) that,  with  $s:=s(n),$ we have $s(n+1)\in \{s,s+1\}.$

Let first  $s(n+1)=s.$ Then $V(n+1)=V(n)$ means $Q(s, n+1)R(s,n+1)=Q(s,n)R(s,n).$ 
Replacing in the proof of  A (i) all inequality signs by "=", this means $$V(n+1)=V(n) \iff  q_{n+1}R(s,n)+p_{n+1}= R(s,n).$$ Since the rhs equation holds  if and only $R(s,n)=1,$ we have proved (a) for the case $s(n+1)=s(n).$
 
 \smallskip 
 
Second, if $s(n+1)=s+1,$ then (13) and (14) in the proof of A (ii) say, when replacing again all inequality signs  by "=",  that $V(n+1)=V(n)$ if and only if one of the following conditions hold

\medskip

($\alpha$)  $q_{s(n)}=q_{n+1}$, that is $p_{s(n)}=p_{n+1}$ and thus $r_{s(n)}=r_{n+1},$ 

\medskip\noindent
or else, from (15) and (17), 

\medskip
($\beta$) $R(s, n)=1+r_{s}.$ 

\medskip\noindent

The condition ($\alpha$) is what we called in (13) the "trivial" case. Since the sequence $(p_j)$ is non-increasing it implies $p_{s(n)}=p_{s(n)+1}= \cdots = p_n=p_{n+1}.$

\smallskip
Concerning condition ($\beta$) in the case  $p_{s(n)}\not=p_{s(n)+1},$ with  $R(s,n)=R(s+1,n)+r_{s}$ we see that it can only hold if $R(s+1,n) =1.$ Since we are in the case $s(n+1)=s+1,$ the condition reads $R(s(n+1),n)=1$, and hence Theorem 4.1 is proved.
\qed

\bigskip 

We note that the Condition (ii) in the theorem is very transparent. If the monotone sequence $(p_j)$ is piecewise constant on a stretch beginning at  $k$ say,
it suffices to check whether the length of the stretch has length at least $1/p_k.$ Condition (i) is in general harder to see (or to exclude) but monotonicity makes it again easier.

\medskip
Essentially the same holds for a corresponding
Theorem for monotone non-decreasing $p_j$'s except
that $s(n+1)$ may now be larger that $s(n) +1.$

\medskip
From the preceding two results we get immediately:

\medskip
\noindent{\bf Corollary 4.2} If $(p_j)$ is non-increasing then optimal thresholds and optimal values of different $n$-problems are all  unique if $R(s(n), n) \not= 1$ and $R(s(n+1), n) \not= 1$ and $p_{s(n)}\not=p_{n+1}$ for all $n$.

\subsection{Strict monotonicity}
For a given underlying sequence $(p_j)$ we say that we have a  {\it coincidence} in $n$ and $n+k,$ if and only if  $V(n)=V(n+k).$ 
Using Theorem 4.1. and Corollary 4.2 we can prove many results about coincidences and their frequencies in an underlying sequence. It turns out that most questions which may come to our mind are seemingly easy to answer.  Here is a small collection:

For example, are there monotone sequences with repeated coincidences $V(n)=V(n+1)=V(n+2)=...$? Clearly yes.
Any sequence which is constant from a certain index $j$ onward, will do. 

Or, can we find a sequence where all $V(n)$ are different from a certain index onward? Yes, an easy choice is one where the sequence of partial sums $S_n:=p_1+p_2+\cdots +p_n$ converges. 

Or a third one, can we find for each monotone sequence $(p_j)$ with diverging partial sums a minored sequence with diverging partial sums such  $\{V(n)=V(n+1)~ {\rm infinitely~often}\}$? 
Again the answer is yes, and first-year analysis suffices for the proof. 
Questions of this kind may not be of general interest, however, and thus this should do.  

Nevertheless, one fact is of  interest, namely, monotone sequences without coincidences yield strictly monotone values $V(n).$ This will be used in  Section 5 (see Example 5.4)

\section{Applications} 
{\bf 5.1}~~Our first example shows that understanding monotonicity may override wrong intuitions. Suppose we are offered 
two games: a 4-game with probabilities $.1, .2, .24, .25$ or the 5-game $.1, .2, .24, .25, .251.$ One may feel it is harder to succeed in the 5-game since the expected number of ones is small  in both games  and the last $1$ has more places to "hide" in the 5-problem with $p_5$ being only slightly larger than $p_4.$ Theorem 3.1  tells us without computation that we should choose the 5-problem, however.  (We gain with probability $.4215$, which is around $1.68$ percent better than playing the 4-game optimally.)

\medskip\noindent
{\bf 5.2}~~To give a different kind of example, we note that an interesting event need not be  associated with a value as such. For instance, in so-called compassionate-use treatments, stopping at the last successful treatment in a sequential clinical trial does not distinguish the last fortunate patient but means stopping with the {\it first} patient covering {\it all} successes,  preventing unnecessary treatments thereafter. The odds of a successful treatment must  typically be estimated sequentially (see B. (2018), 3.2 and 3.3). Here the idea is to plug in the estimates into the odds-algorithm. The independence of the $I_k$ is now lost,  and no optimality can be  claimed, but since different patients are independent of each other in their reaction to treatments, this alternative may still be a good approximation for the optimal strategy. Since new patients tend to join the line if treatments seem more and more successful, or former patients may withdraw from the line if the contrary seems the case, it is good to see for patients and doctors that this is fully in line with Theorem 3.1. 

\medskip
\noindent
{\bf 5.3} May a reasoning based on piecewise monotonicity
also be helpful? Yes, but this depends on where monotonicity begins and on the corresponding $s(n).$ For example, look at the interesting group interview  problem of Hsiau and Yang (2000). (Example 3.3 of B. (2000) shows the concise solution with the odds-algorithm). This is a give-and-take problem in the sense that we can (formally) win by interviewing all candidates together but we can hardly expect to make good  interviews. Suppose we reserve 5 days for seeing 15 candidates, say, and begin for external reasons with group sizes 3 each on the 1st and 2nd day. Since $I_1=p_1=1$ and$p_k<1$ for $k\ge 2$ by definition, we cannot arrange for having increasing $p_k$'s. If we want groups of sizes $2,3,4$ in any order on days $3,4,5$ it turns out best we choose from the $6$ remaining possibilities the schedule $(3,3,4,3,2)$ and get the best with probability $0.448\cdots.$

\bigskip\noindent
{\bf 5.4} The last  example, concerning the classical secretary problem, is given in full detail:

 \subsection*{Return to the classical secretary problem}
 Grau Ribas  (private communication (2018)) instigated by Bayon et al.\,(2018)) asked whether the optimal win probability in the CSP with $n$ candidates is {\it strictly} decreasing  for $n\ge 3$. Neither Grau Ribas, nor the author, nor peers we asked found the statement in the literature. (A reader knowing a source is kindly requested to inform the author.) 
 
 \medskip
 The answer, given by Theorem 5.1 given below is affirmative, and says more. 
 
 Let  $I_k={\bf 1}\{\rm kth~candidate~has~relative ~rank~ 1\}, k=1,2, \cdots.$ It is well known that the $I_k$ are independent with $p_k=P(I_k=1)=1/k.$ The CSP
for the $n$-problem is the problem to stop online with maximum probability on rank 1 of $n$ uniquely rankable candidates, i.e. on the {\it last} indicator $I_k=1$
for $k\le n.$ Note that $V(2)=1/2=V(3)$ so that we confine our interest to $n\ge 3.$

\bigskip

\noindent {\bf Theorem 5.1}: In the classical secretary problem with $n$ candidates:

\medskip
(i) The optimal win probability $V(n)$ is strictly decreasing in $n$ for $n\ge 3.$

\smallskip
(ii) The optimal thresholds $s(n)$ are all unique for $n\ge 3$.

\bigskip \noindent{\bf Proof}: 
 Let $V(n)$ be  the optimal win probability for $n$ candidates. 
 We have $p_n=1/n, q_n=1-p_n=(n-1)/n$ and $r_n=p_n/q_n=1/(n-1).$  Since $(p_n)$ is decreasing Theorem 3.1 (A) implies that $V(n)$ is non-increasing. It follows that $V(n)$ is strictly decreasing for $n\ge 3$ if and only if $V(n)\not=V(n+1)$ for $n\ge 3.$ With $(p_j)$ being strictly decreasing Corollary 4.2 says then that (i) and (ii) hold both at the same time if  $R(s(n),n)\not=1$ and 
 $R(s(n+1),n)\not=1$ for all $n\ge 3.$ 
 
\smallskip\noindent
Using from above $p_n, q_n,$ and $r_n$ in  $ R(k,n)$ we get
\begin{align} R(k,n) =\sum_{j=k}^n \frac{1}{j-1}=\sum_{j=k-1}^{n-1} j^{-1}=H(n-1)-H(k-2),
\end{align} where $H(n)$ denotes the $n$th partial sum of the harmonic series $H(n)=1+2^{-1}+3^{-1}+ \cdots+n^{-1}.$  If we can show \begin{align} \forall \, n\ge 3  {~\rm and} ~1\le k <n-1: H:= H(n)-H(k) \not \in \N {\rm~for~} n \ge 3,\end{align}then  clearly  $R(k,n)\not=1$ for $n\ge 3.$ 

 Indeed, $H$ cannot be integer, and our proof, or a similar one, may be known. It suffices to study the case that $H$ has at least two summands, since otherwise $R(k,n)<1.$
Look at that denominator $j$ with $k+1 \le j \le n$ and $k+1<n$ of the sum \begin{align}H:=H(n)-H(k)=\frac{1} {k+1}+\frac{1}{k+2}+\cdots+\frac{1}{j}+\cdots+\frac{1}{n} \end{align}the prime factorization of which contains the {\it highest} power of $2,$ that is $j=j_0 \,2^\ell$, say.
Since $H$ has at least two summands,  at least one $j$ must be even, and hence $\ell \ge 1.$ Also, note that  $j_0$ is then odd by definition. 

Now, if $H$ in (24) were integer, then $H\times 2^{\ell-1}$ would also be integer. However, if we bring all summands
of $H\times 2^{\ell-1}$ on a least common denominator $D$, then $D$ must be  even since $j$ is  then replaced by $2  j_0.$ The summands in the corresponding numerator are then all even except exactly one, so that the numerator is odd. This is a contradiction, however, since the ratio of an odd number and an even number cannot be integer, and
This proves (i) and (ii) of Theorem 5.1. at the same time.  \qed 

\bigskip
\noindent {\bf 5.5 Prophet value comparisons}

\smallskip
 It is  easy to find for the CSP alternative proofs for the first part showing that $V(n)$ is non-increasing. The following prophet argument in order to compare $V(n)$ and $V(n+1)$ leads to a proof which
may be the among the most elegant ones.  
   
   \smallskip
 Suppose a prophet knows the position of rank $n+1$, the worst candidate. Let $V_P(n+1)$ be the prophet's value. Clearly $V_P(n+1)\ge V(n+1).$ Knowing the position  of rank $(n+1)$ the prophet will ignore it and thus face an equivalent $n$-problem, since the relative ranks on the other $n$ positions do not change. Hence $V(n)=:V_P(n+1)\ge V(n+1).$\qed

\medskip\noindent A similar prophet argument was already used to show that the value $v_n$ of the celebrated Robbins' Problem for $n$ observations  is non-decreasing in $n$ (see B. and Ferguson (1993), or, independently, Assaf and Samuel-Cahn (1996)). 
For a more general formulation  see the notion of a {\it half-prophet} in B. and Ferguson (1996). 

The preceding argument  was easy due to the nice structure of the CSP model where all positions of ranks are equally likely.  In more general settings such prophet tricks are usually more involved. Note that Theorem 3.1 stays useful, since, whatever index $k,$ fixed or random, is singled
out for a prophet's value comparison, cancelling $p_k$ in the underlying sequence does not affect any monotonicity property.  

\bigskip
 \noindent{\bf Remark 5.1} If the number of candidates $n$ in the secretary problem is a random variable $N$ then the optimal strategy is in general no longer a simple threshold strategy. As shown in Presman and Sonin (1972) (the fathers of the secretary problem with unknown $n$) the stopping region may, depending on the law $P(N=n)_{n=1,2, \cdots}$ split into {\it stopping islands}, and they found the corresponding monotonicity criteria.

\medskip\noindent We mention that, independently, the  unified approach for unknown $N$ in continuous time
(B. (1984)) may be seen as an interesting alternative model from the point of view of applicability in real-world problems.

\medskip\bigskip\noindent
{\bf References}

\medskip 
D. Assaf and E. Samuel-Cahn, (1996) {\it The secretary problem: minimizing the expected rank with\,i.i.d. random variables}, Adv.\,Appl. Prob., Vol.\,28, No.\,3,  828-852.

\medskip  
L. Bayon,  P. Fortuny, J. Grau, A. Oller-Marcen, and M. Ruis (2018), {\it The best-or-worst and the post-doc problems}, J. Combin. Optim., Vol. 35, No. 3, 703-723.

\medskip
F.T. Bruss (1984), {\it A Unified Approach to a Class of Best Choice Problems with an Unknown Number of Options}, ~Ann. Probab., Vol. 12, No. 3, 882-889.

\medskip
F.T. Bruss (2000), {\it Sum the Odds to One and Stop}, ~Ann. Probab., Vol. 28, No. 3, 1384 -1391.

\medskip
F.T. Bruss (2003), {\it  A Note on Bounds for the Odds-theorem of Optimal Stopping}, ~Ann. Probab.,
Vol. 31, No. 4, 1859-1861.

\medskip
F.T. Bruss (2018), {\it A mathematical approach to comply with ethical constraints in compassionate use treatments}, Math. Scientist, Vol. 43, No. 1, 10-22.

\medskip
F.T. Bruss and  T.S. Ferguson (1993), {\it Minimizing the expected rank with full information},
J. Appl. Prob. ,Vol. 30 , No. 3, 616-626. 

\medskip
F.T. Bruss and  T.S. Ferguson (1996), {\it Half-prophets and Robbins' Problem of Minimizing the Expected Rank}, Springer Lecture Notes in Statistics, No. 114, Vol. I:  Applied Probability, 1-17.

\medskip
R. Dendievel 10 (2013), {\it New Developments of the Odds-Theorem}, Math. Scientist, Vol. 38, 111-123.

\medskip
T.S. Ferguson (2016), {\it The Sum-the-Odds Theorem with Application to a Stopping Game of Sakaguchi}, Math. Applicanda, Vol. 44 (1), 45-61.

\medskip
A.V. Gnedin (1994), {\it On a best-choice problem by dependent criteria}, J. Appl. Prob., Vol. 31, No. 1, 221-234.

\medskip
A. Goldenshluger, Y. Malinovsky, and A. Zeevi (2019), {\it A Unified Approach for Solving Sequential Selection Problems}, arXiv:1901.04183v3, \,Jan 2019.

\medskip
J.M. Grau Ribas (2018), {\it A New proof and Extension of the Odds-Theorem},\\arxiv:1812.09255v1, \,Dec. 2018.

\medskip
S.-R. Hsiau and J.-R. Yang (2000), {\it A natural variation of the standard secretary problem}, Statistica Sinica, Vol. 10, No. 2, 639-646.

\medskip
T. Matsui and K. Ano, (2014), {\it A Note on a Lower Bound for the Multiplicative Odds-Theorem of Optimal Stopping}, J. Appl. Prob., Vol 51, No. 3, 885-889. 

\medskip
T. Matsui and K. Ano, (2016)
{\it Lower Bounds for Bruss’ Odds Problem with Multiple Stoppings}, Math. of Oper. Res., Vol. 41, No. 2, 700-714.

\medskip
K. Szajowski (2007), {\it A game version of the Cowan-Zabczyk-Bruss problem}, Stat. and Probab. Letters, Vol. 77, No. 17, 1683-1689.

\medskip
E.L. Presman and I.N. Sonin (1972), {\it The best choice problem for a random number of objects}, Theory Probab. Applic., Vol. 17, 657-668.

\medskip
M. Tamaki (2010), {\it Sum the multiplicative odds to one and stop,} J. Appl. Probab., Vol. 47, No. 3, 761–777.

\medskip
M. Tamaki (2011), {\it Maximizing the probability of stopping on any of the last m successes in independent Bernoulli trials with random horizon}, Adv. Appl. Prob., Vol. 43, No. 3, 760-781.

\medskip
\centerline{---}

\medskip\noindent{\bf Author's address}:\\ F.\,Thomas Bruss, \\Universit\'e Libre de Bruxelles,\\Facult\'e des Sciences, CP 210, \\B-1050 Brussels, Belgium;   (tbruss@ulb.ac.be)

\end{document}